\documentclass[12pt]{article}
\usepackage{amsmath}
\usepackage{amssymb}
\usepackage{amsthm}
\usepackage{amsfonts}
\usepackage{color}
\definecolor{labelkey}{rgb}{1,0,0}
\usepackage{graphicx}
\numberwithin{equation}{subsection}
\newtheorem{theorem}[equation]{Theorem}
\newtheorem{prop}[equation]{Proposition}
\newtheorem{lemma}[equation]{Lemma}
\newtheorem{cor}[equation]{Corollary}

\newtheorem{conj}[equation]{Conjecture}

\theoremstyle{definition}

\theoremstyle{remark}
\newtheorem{remark}[equation]{Remark}

\newcommand{\cx}{{\bf C}}

\newcommand{\RR}{{\Bbb R}}

\newcommand{\cP}{{\cal{P}}}

\newcommand{\SSS}{{\Bbb S}}
\DeclareMathOperator{\Id}{Id}

\newcommand{\Lap}{{\bf \Delta}}

\newcommand{\C}{{\mathbb C}}
\DeclareMathOperator{\Fix}{Fix}
\DeclareMathOperator{\supp}{supp}
\DeclareMathOperator{\grad}{\nabla}
\newcommand{\Tau}{\mathcal{T}}
\newcommand{\RE}{\mathbb{R}}

\newcommand{\area}{{\rm Area}}

\begin{document}
\title{How large can the first eigenvalue be on a surface of genus two?}
\author{
Dmitry Jakobson%
\thanks{Department of Mathematics and Statistics, 
McGill University,
805 Sherbrooke Str. West, Montreal, QC H3A 2K6, Canada; 
e-mail \texttt{jakobson@math.mcgill.ca}}
\and
Michael Levitin%
\thanks{Department of Mathematics, 
Heriot-Watt University,
United Kingdom; 
e-mail  \texttt{M.Levitin@ma.hw.ac.uk}}
\and
Nikolai Nadirashvili%
\thanks{Laboratoire d'Analyse, Topologie, Probabilit\'{e}s UMR 6632,
Centre de Math\'{e}matiques et Informatique, 
Universit\'{e} de Provence, 
39 rue F. Joliot-Curie, 
13453 Marseille Cedex 13,
France;  
e-mail  \texttt{nicolas@cmi.univ-mrs.fr}}
\and
Nilima Nigam%
\thanks{Department of Mathematics and Statistics, 
McGill University,
805 Sherbrooke Str. West, Montreal, QC H3A 2K6, Canada; 
e-mail \texttt{nigam@math.mcgill.ca}}
\and
Iosif Polterovich%
\thanks{D\'{e}partement de math\'{e}matiques et de statistique,
Universit\'e de Montr\'eal, 
CP 6128 succ Centre-Ville,
Montr\'eal, QC H3C 3J7, Canada; 
e-mail \texttt{iossif@dms.umontreal.ca}}
}

\date{11 September 2005}

\maketitle

\begin{abstract}
Sharp upper bounds for the first eigenvalue of the Laplacian on a surface of a fixed area
are known only in genera zero and one. We investigate the genus two case and conjecture that
the first eigenvalue is maximized on a singular surface which is realized as a  double branched covering 
over a sphere. The six ramification points are chosen in such a way that this surface has a complex structure of the Bolza surface. We prove that our conjecture follows from a lower bound on the first eigenvalue
of a certain mixed Dirichlet-Neumann boundary value problem on a half-disk. The latter can be studied
numerically, and we present conclusive evidence supporting the conjecture. 
\end{abstract}

\


\

\section{Introduction and main results}\label{sec:intro}
\subsection{Upper bounds on the first eigenvalue}
Let $M$ be a closed  surface of genus $\gamma$ and let
$g$ be the Riemannian metric on $M$. Denote by
$\Delta$ the
Laplace-Beltrami operator on $M$, and by $\lambda_1$
the smallest
positive eigenvalue of the Laplacian. Let the area
$\area(M)$ be fixed. 
How large can $\lambda_1$ be  on such a surface?

Sharp bounds for the first eigenvalue are known only
for the sphere (\cite{H}, see also \cite{SY}),
the projective plane (\cite{LY}), the torus
(\cite{B}, \cite{N1}), and the Klein bottle (\cite{JNP1}, 
\cite{EGJ}).   The present paper is concerned with the
surface of genus $2$.

Let $M$ be orientable and let $\Pi: M \to \mathbb{S}^2$ be a 
non-constant holomorphic map (or, conformal
branched covering) of degree $d$. It was proved in
\cite{YY} that
\begin{equation}
\label{Yang}
\lambda_1 \area(M) \le 8\pi d\,.
\end{equation}
Any Riemann surface of genus $\gamma$ can
be represented as a branched cover over $\mathbb{S}^2$
of degree
$\displaystyle d=\left[\frac{\gamma+3}{2}\right]$,
where $[\cdot]$ denotes the
integer  part (see \cite{Gun}, \cite{GH}). Therefore, 
\begin{equation}\label{upperbd}
\lambda_1\area(M) \le 8\pi
\left[\frac{\gamma+3}{2}\right].
\end{equation}
In general, \eqref{upperbd} is not sharp, for example
for $\gamma=1$ (\cite{N1}).

Let  $M=\cP$ be a surface of genus $\gamma=2$.
Then \eqref{upperbd} implies
\begin{equation}
\label{sharp}
\lambda_1\area(\cP) \le 16\pi.
\end{equation}
The aim of this paper is to show, using a mixture of
analytic and numerical tools, 
that (\ref{sharp}) is sharp. Main results of this paper were announced (without proofs)
in \cite[section 4]{JLNP}. 

\subsection{The Bolza surface}
Let $\Pi:\cP \to \SSS^2$ be a branched covering of degree $d=2$.
The Riemann-Hurwitz formula (see \cite{GH}) implies that this cover is
ramified at $6$ points.  We choose these points to be the
intersections of the round sphere $\SSS^2$ centered at the origin
with the coordinate axes in $\RR^3$. The surface $\cP$ can
be realized as
$$
\left\{(z,w)\in\cx^2:w^2=F(z):=z\frac{(z-1)(z-i)}{(z+1)(z+i)}\right\}\,.
$$
This surface has the conformal structure of the Bolza surface. It has an octahedral group of holomorphic 
automorphisms and its symmetry group is the largest among  surfaces of genus two \cite{Igusa, KW}. 
Interestingly enough, the Bolza surface appears in some other extremal problems, in particular for systoles 
(see \cite{Katz}). 

To simplify calculations it is convenient to rotate the equatorial plane by $\pi/4$.  
The equation of $\cP$ becomes
\begin{equation}\label{proj:def}
\cP:=\left\{(z,w)\in\cx^2:w^2=F(z):=z \frac{(z-e^{\pi
i/4})(z-e^{3\pi i/4})}{(z+e^{\pi i/4})(z+e^{3\pi i/4})}\right\}\,.
\end{equation}
The projection $\Pi$ is defined by $\Pi:(z,w)\to z$. The
set of ramification points in the complex $z$ plane is
$R:=\{0,\infty,\pm e^{\pi i/4},\pm e^{3\pi i/4}\}$.  The spherical
and complex models are related by the stereographic projection;
the induced metric in the complex plane (which we assume coincides
with the equatorial plane of $\SSS^2$) is
\begin{equation}\label{eq:metr}
4dzd\bar{z}/(1+|z|^2)^2.
\end{equation}

Let $g_0$ be the
metric on $\cP$  which is the pullback of the round metric \eqref{eq:metr}
on $\SSS^2$. One can see that the metric $g_0$ has conical
singularities at the points of ramification. It has curvature $+1$
everywhere except the branching points. Because of the presence of singularities we have to
specify what we mean by the first positive eigenvalue of the Laplacian on $(\cP, g_0)$.
We set
$$
\lambda_1(\cP, g_0):=\inf_{u\in H^1_0(\cP, g_0)\,,\ u\ne 0\,,\ \langle u,1\rangle=0}\frac{\|\grad u\|^2}{\|u\|^2}\,,
$$
where the scalar product  $\langle\cdot,\cdot\rangle$ and the norm $\|\cdot\|$ are taken in the space $L_2(\cP, g_0)$.
The Sobolev space $H^1_0(\cP, g_0)$ of functions supported away from the singularities is obtained by the closure 
of $C_0^\infty(\cP, g_0):=\{v\in C^\infty(\cP, g_0)\mid \overline{\Pi\supp v}\,\cap\, R=\emptyset\}$ with respect 
to the norm $\|\grad v\|^2+\|v\|^2$. 

\subsection{Main results}\label{sec:mainres}

We start with the following

\begin{conj}\label{conj1}
The equality in \eqref{sharp} is attained for the metric $g_0$
on $\cP$, i.e.
$$
\lambda_1(\cP, g_0)\area(\cP, g_0)=16\pi.
$$
\end{conj}
Since $(\cP, g_0)$ is a double cover of the standard $\SSS^2$, we
have
$$\area(\cP, g_0)=2\area(\SSS^2)=8\pi.$$
Therefore, in order to prove Conjecture~\ref{conj1} it suffices to
show that
\begin{equation}\label{lambda}
\lambda_1 (\cP, g_0)=\lambda_1(\SSS^2)=2.
\end{equation}

Unfortunately, we are unable to prove \eqref{lambda}, and therefore establish
Conjecture~\ref{conj1}. We can however reduce the conjecture to the following
spectral problem on a quarter-sphere $Q\subset\SSS^2$ that can be treated using numerical methods.
Namely, let,
in usual spherical coordinates $(\phi,\theta)$,
$$
Q = \{(\phi,\theta):\ 0<\phi<\pi/2, 0<\theta<\pi\}\,.
$$
We split the boundary $\partial Q$ into two parts:
$\partial Q = \overline{\partial_1 Q \sqcup \partial_2 Q}$, where
$$
\partial_1 Q = \{(0,\theta): |\theta-\pi/2|<\pi/4\}\cup
\{(\pi/2,\theta):\ 0<\theta<\pi/2\}\,,
$$
$$
\partial_2 Q = \{(0,\theta):\ |\theta-\pi/2|>\pi/4\}\cup
\{(\pi/2,\theta):\pi/2<\theta<\pi\}\,,
$$
and consider the spectral boundary value problem for the
Laplace-Beltrami operator on $Q$:
\begin{equation}\label{eq:sphQ}
-\Delta u = \Lambda u\quad\text{on }Q\,,\quad u|_{\partial_1 Q} = 0\,,
\quad(\partial u/\partial n)|_{\partial_2 Q} = 0\,.
\end{equation}

Let $\Lambda_1$ denote the first eigenvalue of the problem \eqref{eq:sphQ}
(which we understand as usual in the variational sense).

\begin{conj}\label{conj2}
$$
\Lambda_1\ge 2.
$$
\end{conj}

Our main result is
\begin{theorem}\label{conj2to1}
Conjecture \ref{conj2} implies Conjecture \ref{conj1}.
\end{theorem}
Theorem \ref{conj2to1} is proved in sections 2 and 3.

Extensive numerical
calculations (see Section \ref{sec:num}) show that $\Lambda_1 \gtrsim 2.27$ which implies Conjecture \ref{conj1}.
The best lower bound we are able to prove is just $\Lambda_1>0.75$, which follows from Dirichlet-Neumann bracketing
(replace the Dirichlet condition by the Neumann one on the arc $(0,\theta)$, $\pi/4<\theta<3\pi/4)$,   see \cite{Keller}).

Finally, we note that the spectral problem \eqref{eq:sphQ} easily
reduces via the stereographic projection to the following mixed
Dirichlet-Neumann problem on a half-disk  $D:=\{(r, \psi)\in \RR^2: r<1,
0<\psi<\pi\}$ (here $(r,\psi)$ are usual planar polar coordinates):
\begin{equation}\label{eq:polD}
-\Delta v = \frac{4\Lambda}{(1+r^2)^2} v\quad\text{on }D\,,\quad
v|_{\partial_1 D} = 0\,,
\quad(\partial v/\partial n)|_{\partial_2 D} = 0\,.
\end{equation}
Here $\partial_1 D:=\{(r,0): r\in(0,1)\}\cup\{(1,\psi): |\psi-\pi/2|<\pi/4\}$ and
$\partial_2 D:=\{(r,\pi): r\in(0,1)\}\cup\{(1,\psi): |\psi-\pi/2|>\pi/4\}$. 

Problems \eqref{eq:sphQ} and \eqref{eq:polD} are quite
remarkable in their own right --- each of them is an example of a mixed
Dirichlet-Neumann boundary value problem whose spectrum is invariant under a
swap of Dirichlet and Neumann boundary conditions.  Namely, the spectrum of \eqref{eq:polD}
coincides with the spectrum of 
\begin{equation}\label{eq:polD1}
-\Delta v = \frac{4\Lambda}{(1+r^2)^2} v\quad\text{on }D\,,\quad
v|_{\partial_2 D} = 0\,,
\quad(\partial v/\partial n)|_{\partial_1 D} = 0\,.
\end{equation}
We refer to \cite{JLNP} for a further discussion on Dirichlet-Neumann swap isospectrality.

\begin{figure}[!hbt]
\begin{center}
\includegraphics[width=0.45\textwidth]{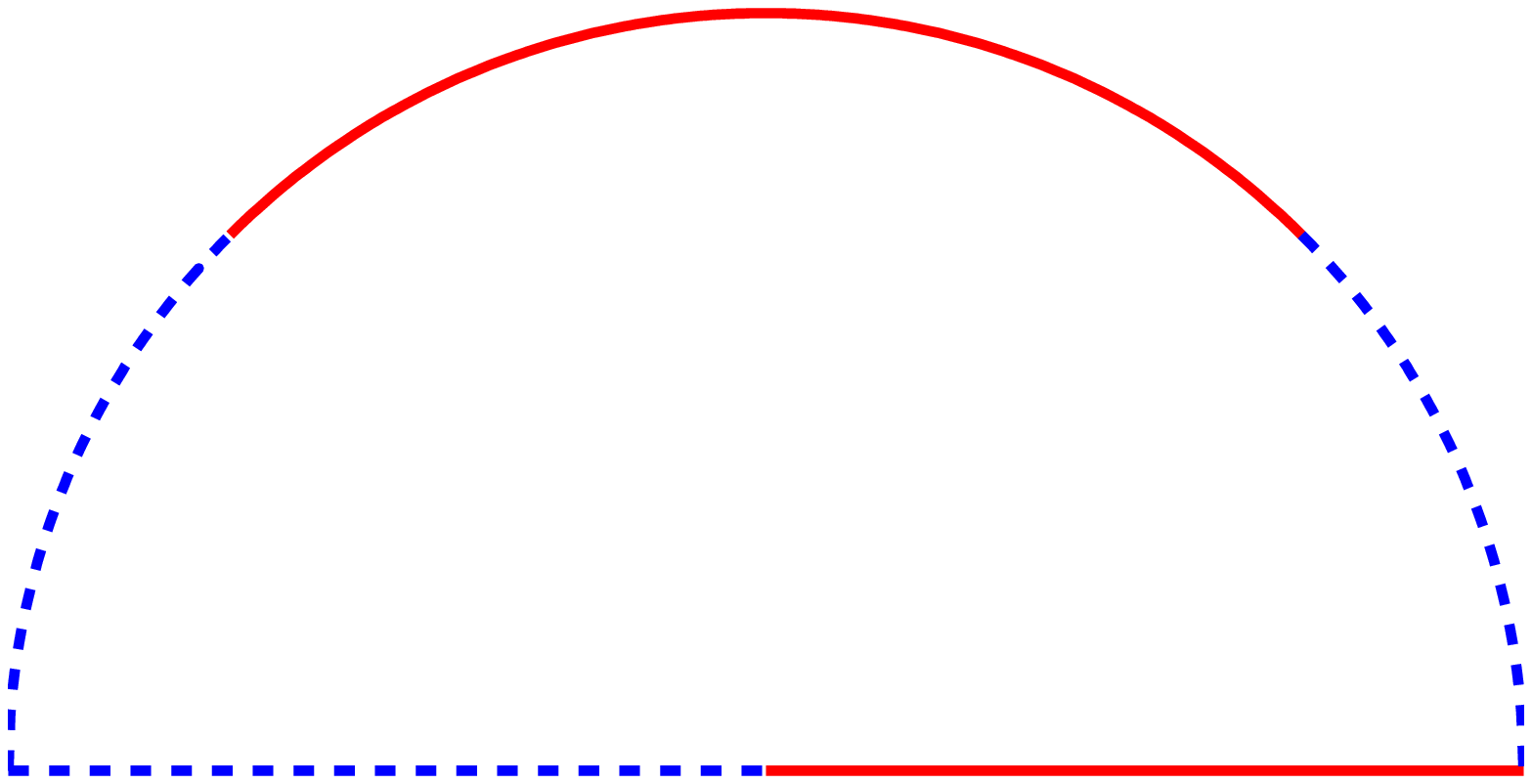}\hfill
\includegraphics[width=0.45\textwidth]{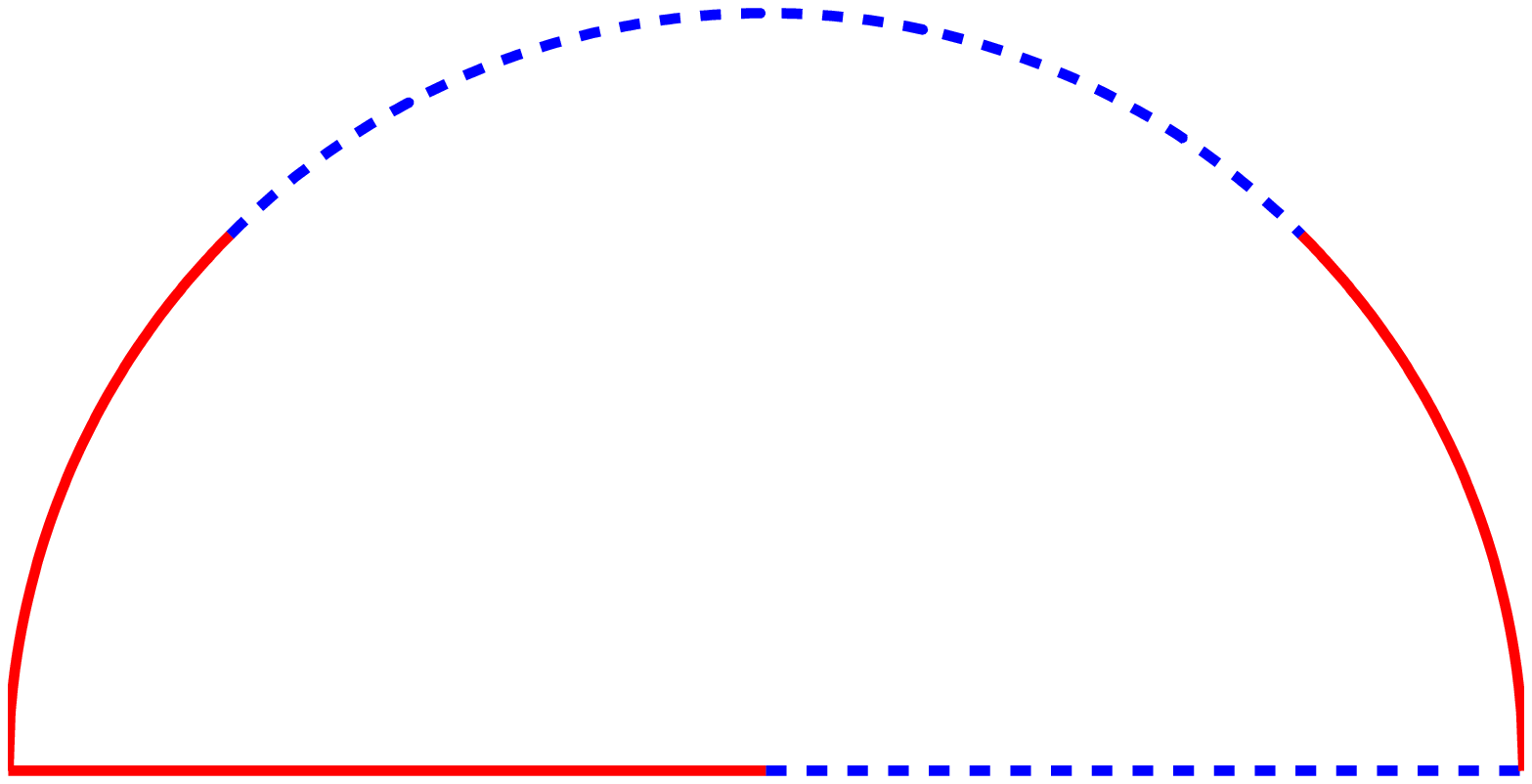}
\caption{Geometry of boundary value problems \eqref{eq:polD} (left) and \eqref{eq:polD1} (right). 
Here and further on, the solid red line denotes Dirichlet boundary condition and the dashed blue line --- 
the Neumann one.}
\label{fig:halfdisks}
\end{center}
\end{figure}

\begin{remark}
One can check that a surface with a finite number of conical singularities can be approximated by a sequence of smooth surfaces 
of the same genus and area in such a way that the corresponding sequence of the first non-zero eigenvalues converges to $\lambda_1$
on the original surface. Thus, Conjecture~\ref{conj1} means that \eqref{sharp} is sharp in the class of smooth metrics, 
although the equality is not necessarily attained. For a general result about the convergence 
of the whole spectrum see \cite{Ro}.
\end{remark}

\section{Symmetries}\label{sec:symm}
\subsection{Hyperelliptic involution}
Let $T: \cP \to \cP$, $T^2 =\Id$ be a map that intertwines the
preimages of points of $\SSS^2$ under a two-sheeted covering
$\Pi:\cP \to \SSS^2$. Clearly,  the Laplace operator $\Delta$
commutes with $T$.

By the spectral theorem, we can consider separately the restrictions of the
Laplacian onto the spaces of functions which are either even or odd with
respect to $T$. The even functions on $\cP$ can be identified with the
functions on $\SSS^2$. Therefore, as $\lambda_1(\SSS^2)=2$, we have
$\lambda_1(\cP) \le 2$, and the equality in  \eqref{lambda} will be achieved if
and only if the first eigenvalue  $\lambda_1^{\text{odd}}$ of the Laplacian
acting on the odd subspace  satisfies $\lambda_1^{\text{odd}} \ge 2$.

\subsection{Isometries of $\cP$}
Consider the following isometries of $\SSS^2$ (as usual, we identify $\SSS^2$
and $\C$ by stereographic projection):
\begin{equation}\label{symm:base}
\begin{aligned}
\sigma_1 &:z \mapsto \bar{z}\quad\text{or}\quad(\chi,\eta,\xi)\mapsto(\chi,-\eta,\xi),\\
\sigma_2 &:z \mapsto -\bar{z}\quad\text{or}\quad(\chi,\eta,\xi)\mapsto(-\chi,\eta,\xi),\\
\sigma_3 &:z \mapsto 1/\bar{z}\quad\text{or}\quad(\chi,\eta,\xi)\mapsto(\chi,\eta,-\xi).
\end{aligned}
\end{equation}
Here $z=x+iy$ is a point in the equatorial plane upon which a point
$(\chi,\eta,\xi)\in\SSS^2$ is stereographically projected.

The hyperelliptic involution $T$ is given by $T:(z,w)\to (z,-w)$.
For $1\leq j\leq 3$,  a symmetry $\sigma_j$ of $\SSS^2$ has two corresponding
symmetries $s_j$ and $T\circ s_j$ satisfying
\begin{equation}\label{hg:rel}
\Pi\circ s_{j}=\Pi\circ T\circ s_{j}=\sigma_j\circ\Pi\,.
\end{equation}
Those symmetries, with account of \eqref{proj:def} are given by the explicit formulae
\begin{equation}\label{symm:cover}
\begin{aligned}
s_1:(z,w)&\mapsto (\bar{z},\bar{z}/\bar{w}),\\
s_2:(z,w)&\mapsto (-\bar{z},i\bar{w}),\\
s_3:(z,w)&\mapsto (1/\bar{z},\bar{w}/\bar{z}).
\end{aligned}
\end{equation}
As an illustration, we demonstrate how the last of these formulae is obtained: if
$w^2 = F(z)$, then by \eqref{proj:def},
$$
F\left(\frac{1}{\bar{z}}\right) =  \frac{1}{\bar{z}}\frac{(1/\bar{z}-e^{\pi
i/4})(1/\bar{z}-e^{3\pi i/4})}{(1/\bar{z}+e^{\pi i/4})(1/\bar{z}+e^{3\pi i/4})}=
\frac{\overline{F(z)}}{\bar{z}^2}\,,
$$
thus giving the expression for $s_3$.

It easily seen that all $s_j$ commute with $T$ and satisfy
\begin{equation}\label{ident:1}
\begin{aligned}
s_j^2&=\Id,\ j=1,2,3;\\
s_1s_3&=s_3s_1,\; s_2s_3=s_3s_2;\\
s_2s_1&=Ts_1s_2.
\end{aligned}
\end{equation}

\begin{remark}
In the proof of Theorem \ref{conj2to1} we will use only the symmetries $s_1$ and $s_3$.  Calculations for $s_2$ 
are presented for
the sake of completeness (see Remark \ref{rems2}).
\end{remark}
\subsection{Fixed point sets of isometries}\label{sec:fixedpt}
Let $\Fix(S)$ denote a fixed point set of a mapping $S$.
As easily seen from \eqref{symm:base},
the sets $\Fix(\sigma_j)$, for $j=1,2,3$, lie in the union  of the coordinate lines and a unit circle of $\C$, and
we introduce the following notation for future reference.  The coordinate lines are divided into two rays
each by the ramification point $r_0:=0$, and we denote
$$
a_1:=\{z=t,t>0\}\,,\quad a_2:=\{z=it,t>0\}\,,
$$
$$
a_3:=\{z=t,t<0\}\,,\quad a_4:=\{z=it,t<0\}\,.
$$
The circle is divided into four arcs by the ramification points
$r_1:=e^{-\pi i/4}$, $r_2:=e^{\pi i/4}$, $r_3:=e^{3\pi i/4}$, and $r_4:=e^{-3\pi i/4}$, and we denote
the arcs by
$$
a_{k+4}:=\{z=e^{t\pi i/4}, t\in(2k-3,2k-1)\}\,,\qquad k=1,2,3,4\,,
$$
so that the arc $a_5$ goes from $r_1$ to $r_2$, the arc $a_6$ goes from $r_2$ to $r_3$, the arc  $a_7$
goes from $r_3$ to $r_4$, and finally $a_8$ goes from $r_4$ to $r_1$.

In this notation, the fixed point sets $\Fix(\sigma_j)$ are written as
\begin{equation}\label{fixset:base}
\Fix(\sigma_1)=a_1\cup a_3\,,\quad
\Fix(\sigma_2)=a_2\cup a_4\,,\quad
\Fix(\sigma_3)=a_5\cup a_6\cup a_7\cup a_8\,.
\end{equation}

\begin{figure}[!hbt]
\begin{center}
\includegraphics[width=0.75\textwidth]{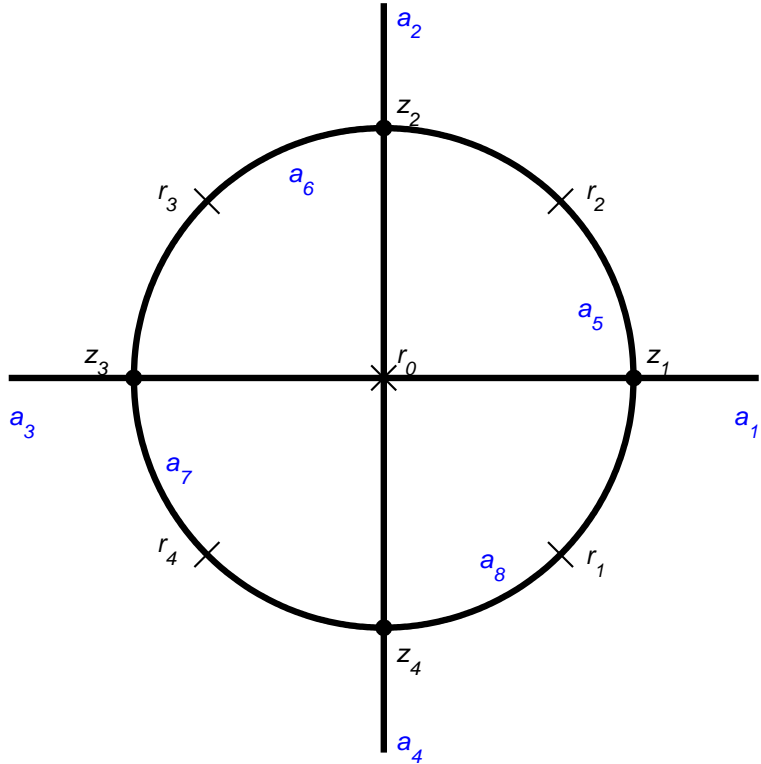}
\caption{Ramification points, rays, arcs and intersections}
\label{fig:points}
\end{center}
\end{figure}

Note that each of the rays $a_j$ ($j=1,2,3,4$) intersects an arc $a_{j+4}$ at a single
point which we
denote $z_j$:
$$
z_1 = 1\,,\quad z_2=i\,,\quad z_3=-1\,,\quad z_4=-i\,,
$$  
see Figure \ref{fig:points}.

\subsection{Fixed point sets of $s_1,s_2,s_3$}\label{subs:fixa1}
Each of the points $z_j$ has exactly two pre-images
$p_j^{(m)}:=(w_j^{(m)}, z_j)\in\Pi^{-1}z_j$, $m=1,2$, where $w_j^{(1,2)}$ are the
solutions of the equation $(w_j)^2=F(z_j)$, with $F$ given in \eqref{proj:def}. These
solutions are easily found from  \eqref{proj:def};
we are of course at liberty to choose which of the two solutions is denoted $w_j^{(1)}$
 and which is denoted
$w_j^{(2)}$. For definiteness we set
\begin{equation}\label{eq:ws}
\begin{split}
w_1^{(1)}=i\,,\qquad& w_1^{(2)}=-i\,;\\
w_2^{(1)}=\frac{1+i}{2+\sqrt{2}}\,,\qquad& w_2^{(2)}=-\frac{1+i}{2+\sqrt{2}}\,;\\
w_3^{(1)}=1\,,\qquad& w_3^{(2)}=-1\,;\\
w_4^{(1)}=\frac{1-i}{2-\sqrt{2}}\,,\qquad& w_4^{(2)}=-\frac{1-i}{2-\sqrt{2}}\,.
\end{split}
\end{equation}

For future use, we need to know the images of points $p_j^{(m)}$ under the symmetries
$s_l$, $l=1,2,3$.
These are easily calculated from \eqref{symm:cover}; it turns out that
$s_lp_j^{(m)}=p_k^{(n)}$ with some indices $k\in\{1,2,3,4\}$, $n\in\{1,2\}$. The
results of the calculations are summarized in the following Table~\ref{tab:ws_and_ss}.

\begin{table}[htdp]
\renewcommand{\tabcolsep}{0.5cm}
\renewcommand{\arraystretch}{1.3}
\begin{center}
\begin{tabular}{|c||c|c|c|}
\hline
$(j,m)$&\multicolumn{3}{c|}{$(k,n)$}\\\cline{2-4}
&$l=1$&$l=2$&$l=3$\\\hline
(1,1)&(1,1)&(3,1)&(1,2)\\
(1,2)&(1,2)&(3,2)&(1,1)\\\hline
(2,1)&(4,1)&(2,1)&(2,1)\\
(2,2)&(4,2)&(2,2)&(2,2)\\\hline
(3,1)&(3,2)&(1,2)&(3,2)\\
(3,2)&(3,1)&(1,1)&(3,1)\\\hline
(4,1)&(2,1)&(4,2)&(4,1)\\
(4,2)&(2,2)&(4,1)&(4,2)\\\hline
\end{tabular}
\caption{The points $p_j^{(m)}$ and their images $p_k^{(n)}=s_l p_j^{(m)}$ under symmetries $s_l$. The table lists the pairs of indices 
$(j,m)$ and the resulting pairs $(k,n)$ for $l=1,2,3$. Note also that $T$ acts by interchanging the second indices $1\leftrightarrow 2$.\label{tab:ws_and_ss}}
\end{center}
\end{table}

For each of the rays or arcs $a_k$, $k=1,\dots,8$ defined in the previous section,
its pre-image $\Pi^{-1}a_k$ has two connected components which we denote
$b_k^{(1)}, b_k^{(2)}$ related by $b_k^{(m)}=Tb_k^{(m')}$, $m, m'\in\{1,2\}$, $m\ne m'$.
Again, the choice of which component we denote by an upper index $(1)$ is up to us and
in order to fix the notation we postulate that $b_k^{(m)}\ni w_{k'}^{(m)}$, $k'=((k-1)\mod 4)+1$, e.g.
$w_{1}^{(1)}\in b_1^{(1)}\cap b_5^{(1)}$, $w_{3}^{(2)}\in b_3^{(2)}\cap b_7^{(2)}$, etc.

We now have at our disposal all the information  we need in order to obtain the
fixed point sets of $s_1, s_2, s_3$. We start with the following two simple Lemmas.

\begin{lemma}\label{lem:ssig1}  $\Pi\Fix(s_j)\subseteq\Fix(\sigma_j)$.
\end{lemma}

\begin{proof} Let $z\in\Pi\Fix(s_j)$. Then there exists $p\in\cP$ such that $\Pi p=z$ and 
$s_j p=p$. Thus $\Pi s_j p =z$ and by \eqref{hg:rel} $\sigma_j \Pi p = \sigma_j z = z$, so 
that $z\in\Fix(\sigma_j)$.
\end{proof}

\begin{lemma}\label{lem:ssig2} Let $a_k\subseteq\Fix(\sigma_j)$. Then, for $m=1,2$,  either
$b_k^{(m)}\subseteq\Fix(s_j)$ or
$b_k^{(m)}\subseteq\Fix(T\circ s_j)$.
\end{lemma}

\begin{proof}
We have $\Pi b_k^{(m)}=a_k$, so that $\sigma_j\Pi b_k^{(m)}=\sigma_j a_k=a_k$, and so
by \eqref{hg:rel}
$\Pi s_j b_k^{(m)}= a_k=\Pi b_k^{(m)} = \Pi T b_k^{(m)}$. The result follows from the
obvious observation:
if $\Pi\alpha=\Pi\beta$, then either $\alpha=\beta$ or $\alpha=T\beta$.
\end{proof}

The lemmas lead  to the following
\begin{prop}\label{prop:fix}
\begin{equation*}
\begin{split}
\Fix(s_1) &= b_1^{(1)}\cup b_1^{(2)}\,,\\
\Fix(Ts_1) &= b_3^{(1)}\cup b_3^{(2)}\,,\\
\Fix(s_2) &= b_2^{(1)}\cup b_2^{(2)}\,,\\
\Fix(Ts_2) &= b_4^{(1)}\cup b_4^{(2)}\,,\\
\Fix(s_3) &= b_6^{(1)}\cup b_6^{(2)}\cup b_8^{(1)}\cup b_8^{(2)}\,,\\
\Fix(Ts_3) &= b_5^{(1)}\cup b_5^{(2)}\cup b_7^{(1)}\cup b_7^{(2)}\,.
\end{split}
\end{equation*}
\end{prop}

\begin{proof} By Lemmas \ref{lem:ssig1} and \ref{lem:ssig2}, for any given $j$ the fixed sets $\Fix(s_j)$ and $\Fix(Ts_j)$ consist only of 
the pre-images of the components $a_k$ of the corresponding fixed sets $\Fix(\sigma_j)$ (given by \eqref{fixset:base}). However we still 
need to describe which component $b_k^{(m)}$, $m=1,2$, lies in 
$\Fix(s_j)$ and which in $\Fix(Ts_j)$. As each component $b_k(m)$ is uniquely determined by the point $w_k^{(m)}$ given by \eqref{eq:ws}, 
it is sufficient just to check in Table~\ref{tab:ws_and_ss} whether  $s_j w_k^{(m)}=w_k^{(m)}$ or $Ts_j w_k^{(m)}=w_k^{(m)}$.

For example, to find $\Fix(s_2)$ we need only to inspect $b_2^{(m)}$ and $b_4^{(m)}$. As, by Table \ref{tab:ws_and_ss},  
$s_2 w_2^{(m)}=w_2^{(m)}$ and 
$Ts_2 w_4^{(m)}=w_4^{(m)}$, we have $\Fix(s_2) = b_2^{(1)}\cup b_2^{(2)}$ and
$\Fix(Ts_2) = b_4^{(1)}\cup b_4^{(2)}$. The rest of Proposition \ref{prop:fix} is obtained 
in the same manner.
\end{proof}


\section{Proof of Theorem \ref{conj2to1}}\label{sec:odd}
We divide the proof of  Theorem \ref{conj2to1} into several steps.
\subsection{Even eigenfunctions with respect to $T$}\label{subs:even}
Consider the subspace $V_+\subset L^2(\cP)$ consisting of all even
eigenfunctions with respect to $T$. Any such eigenfunction has a well-defined
projection on $\SSS^2$.  Therefore, if there exists a first eigenfunction of $\cP$
that belongs to $V_+$, its projection is an eigenfunction on $\SSS^2$ and hence
the corresponding eigenvalue is greater or equal than two (recall that $\lambda_1(\SSS^2)=2$).
Hence, in this case the Conjecture \ref{conj1} is verified.

\subsection{Use of symmetries $s_1,s_3$.}\label{symm:135}
Denote by $G_{13}$ the subgroup of the automorphism group of $\cP$
generated by the symmetries $$\{T,s_1,s_3\}.$$  
It follows from \eqref{ident:1} that $G_{13}$ is commutative. 
Note also that all the elements of $G_{13}$ commute with the Laplacian on $\cP$.
Therefore, we can choose a basis of $L^2(\cP)$
consisting of joint eigenfunctions of all $s\in G_{13}$ and $\Delta$.
Given a joint eigenfunction $f$ of all $s\in G_{13}$, we denote by
$\mu(f,s)$ the corresponding eigenvalue of $s$, i.e.
$$
f(sx)=\mu(f,s)f(x)\,.
$$
Since $s_j^2=T^2=\rm{Id}$ for $j=1,3$, we see that $\mu(f,s)=\pm 1$ for all $s\in G_{13}$.

\subsection{Odd eigenfunctions with respect to $T$}
Consider now the space $V_-\subset L^2(\cP)$ consisting of all 
eigenfunctions of the Laplacian which are {\em
odd} with respect to $T$. 
Let $\phi_1$ be a joint eigenfunction of $\{T, s_1, s_3, \Delta\}$, 
corresponding to the smallest eigenvalue of $\Lap\left|_{V_-}\right.$ 

Now, since $\mu(\phi_1,T)=-1$ and $s_3^2T=T$, we have
$\mu(\phi_1,s_1)\mu(\phi_1,s_1 T)=\mu(\phi_1,T)=-1$, and similarly
$\mu(\phi_1,s_3)\mu(\phi_1,s_3 T)=-1$.

Without loss of generality we may assume that $\mu(\phi_1,s_1)=-1$. 
We recall from section
\ref{sec:fixedpt} that the fixed point set  $\Fix{s_1}$ consists of
the arcs $b_1^{(1)},b_1^{(2)}$.  Thus $\phi_1$ must {\em vanish} on these arcs.

Consider now the symmetries $s_3,s_3T$. We must have one of the following two cases:
\begin{itemize}
\item[i)] 
$\mu(\phi_1,s_3T)=-1$, $\mu(\phi_1,s_3)=1$; 
\item[ii)] 
$\mu(\phi_1,s_3)=-1$, $\mu(\phi_1,s_3T)=1$.
\end{itemize}

Consider first Case i).

\begin{prop}\label{case1}
In Case {\rm i)} the function $\phi_1$ vanishes on the arcs
$$
b_1^{(1)},b_1^{(2)},b_5^{(1)},b_5^{(2)},b_7^{(1)},b_7^{(2)}\,,
$$ 
and its normal derivative $\partial_n\phi_1$ vanishes on the arcs
$$
b_3^{(1)},b_3^{(2)},b_6^{(1)},b_6^{(2)},b_8^{(1)},b_8^{(2)}\,.
$$
\end{prop}
\begin{proof}
By Proposition \ref{prop:fix}, the
fixed-point set of $s_3T$ consists of the arcs
$b_5^1,b_5^2,b_7^1,b_7^2$.  Accordingly, $\phi_1$ vanishes on all
those arcs, as well as on $b_1^1,b_1^2$.
Moreover, $\phi_1$ has $\mu(\phi_1,s_3)=\mu(\phi_1,s_1T)=1$.  It follows that
the normal derivative of $\partial_n\phi_1$ vanishes on the
fixed-point sets of those symmetries.  It remains to apply once more Proposition \ref{prop:fix} in order to complete the proof.
\end{proof}

Consider next Case ii). 
\begin{prop}\label{case2}
In Case {\rm ii)} the function $\phi_1$ vanishes on the arcs
$$
b_3^{(1)},b_3^{(2)},b_6^{(1)},b_6^{(2)},b_8^{(1)},b_8^{(2)}\,,
$$ 
and its normal derivative
$\partial_n\phi_1$ vanishes on the arcs
$$
b_1^{(1)},b_1^{(2)},b_5^{(1)},b_5^{(2)},b_7^{(1)},b_7^{(2)}\,.
$$
\end{prop}
Proposition \ref{case2} is proved in the same way as Proposition \ref{case1}.

\subsection{Final step of the proof}
Since $\phi_1$ is an odd function with respect to the hyperelliptic involution $T$,
its projection upon $\SSS^2$ is not well-defined. However, the projection of $|\phi_1|$
to $\SSS^2$ is well-defined. Denote it by $\psi_1$.  

In Case i), the function $\psi_1$ 
can be chosen as a test function for  the mixed Dirichlet-Neumann boundary 
value problem \eqref{eq:sphQ}.  Assume now Conjecture \ref{conj2} is true and
the first eigenvalue of  \eqref{eq:sphQ} satisfies $\Lambda_1 \ge 2$. Then the Rayleigh
quotient of $\psi_1$ and hence of $\phi_1$ satisfies the same inequality.  But this means that $\psi_1$ cannot be the first 
eigenfunction on $\cP$ since we get a contradiction
with \eqref{sharp}. Therefore, the first eigenfunction of $\cP$ is even with respect
to $T$, and as was shown in section \ref{subs:even} this implies Conjecture \ref{conj1}.

Similarly, in Case ii), the function $\psi_1$ can be chosen as a test function for the mixed
Dirichlet-Neumann boundary volume problem which is obtained from \eqref{eq:sphQ}
by swapping the Dirichlet and the Neumann conditions. However, it was shown in
\cite{JLNP} that this problem is isospectral to \eqref{eq:sphQ}. Therefore, repeating
the same arguments as above we prove that Conjecture \ref{conj1} holds.
This completes the proof of  Theorem \ref{conj2to1}. \qed

\begin{remark} \label{rems2} 
In the proof of Theorem \ref{conj2to1} we have used only the 
symmetries $s_1$ and $s_3$. Alternatively, we could have used $s_2$ and $s_3$.
One can check directly using Proposition \ref{prop:fix} that applying  $s_2$ one obtains
a mixed Dirichlet-Neumann boundary value problem which is equivalent to  \eqref{eq:sphQ} and
hence no additional information about the first eigenfunction is obtained.
\end{remark}

\subsection{A family of extremal surfaces of genus two}
The purpose of this section is to prove the following
\begin{cor}
\label{per}
Conjecture \ref{conj2} implies that there exists a continuous family $\cP_t$ 
of surfaces of genus $2$ such that $\lambda_1 \area(\cP_t)=16\pi$.
\end{cor}
\begin{proof}
Consider the Riemann surface $\cP_t$ defined by the equation
$$
\left\{(z,w):w^2=
\frac{z\left(z-e^{i(\pi/2-t)}\right)\left(z-e^{i(\pi/2+t)}\right)}%
{\left(z-e^{-i(\pi/2-t)}\right)\left(z-e^{-i(\pi/2+t)}\right)}\right\}
$$
where $t\in(0,\pi/2)$. Note that $\cP_{\pi/4}=\cP$.  It is easy to see that  for any $t$,
$\cP_t$ is symmetric with respect to  $s_1$ and $s_3$. Arguing in the same way
as in the proof of Theorem \ref{conj2to1} and using a stereographic projection, we reduce the problem 
on $\cP_t$ to the following two mixed Dirichlet-Neumann boundary value problems on the half-disk $D$:
\begin{equation}\label{eq:polDt}
-\Delta v = \frac{4\Lambda}{(1+r^2)^2} v\quad\text{on }D\,,\quad
v|_{\partial_1(t)} = 0\,,
\quad(\partial v/\partial n)|_{\partial_2(t)} = 0\,.
\end{equation}
and
\begin{equation}\label{eq:polDt1}
-\Delta v = \frac{4\Lambda}{(1+r^2)^2} v\quad\text{on }D\,,\quad
v|_{\partial_2(t)} = 0\,,
\quad(\partial v/\partial n)|_{\partial_1(t)} = 0\,.
\end{equation}

Here $\partial_1(t):=\{(r,0): r\in(0,1)\}\cup\{(1,\psi): |\psi-\pi/2|<t\}$ and
$\partial(t) D:=\{(r,\pi): r\in(0,1)\}\cup\{(1,\psi): \pi/2>|\psi-\pi/2|>t\}$.

We remark that for $t \neq \pi/4$ these two problems are not isospectral. Using
Dirichlet-Neumann bracketing it is easy to see that \eqref{eq:polDt} has a smaller
first eigenvalue than \eqref{eq:polDt1} if $t<\pi/4$ and a larger one if $t>\pi/4$. 
Denote the minimal first eigenvalue of the two problems by $\Lambda_1(t)$.
According to Conjecture 2 and numerical calculations, $\Lambda_1(\pi/4)>2$. 
Since the first eigenvalues of both problems depend continuously and monotonically on parameter $t$, 
and since $\Lambda_1(0)=\Lambda_1(\pi/2)=0.75$  (see section \ref{sec:mainres}), 
there exist numbers $t_1^*\in(0,\pi/4)$ and $t_2^*\in(\pi/4,\pi/2)$ such that $\Lambda_1(t_1^*)=\Lambda_1(t_2^*)=2$ and so
$\Lambda_1(t)\ge 2$ for $t\in [t_1^*, t_2^*]$. Arguing is above,
we deduce that for all surfaces $\cP_t$ corresponding to these values of $t$, estimate
\eqref{sharp} is sharp. This completes the proof of the theorem. 
\end{proof}

Corollary \ref{per} implies that $16\pi$ is a {\em degenerate} maximum for
$\lambda_1\area(M)$ for surfaces of genus two. This is not the case for surfaces of lower genus
on which the metric maximizing the first eigenvalue is unique. Note also that the extremal metrics
in genera zero and one are analytic, while the surfaces $\cP_t$ have singular points.

\section{Numerical investigations}\label{sec:num} 
\subsection{Basics of the Finite Element Method}

In this section we describe the numerical experiments used to estimate the first eigenvalue of \eqref{eq:polD1}.

We define the space $\mathcal{H}$ as the closure of $\{ v\in C^\infty(D) \vert \overline{\supp v}\cap\partial_2 D=\emptyset\}$ with respect 
to the $H^1(D)$ norm. 

The variational setting for the eigenvalue problem is to find the smallest eigenvalue $\lambda\in \mathbb{R}$,  and the associated 
eigenvector $v\in \mathcal{H}$ such that for all $u\in \mathcal{H}$, 
\begin{equation}\label{eq:varP} 
\int_D\nabla v\cdot \nabla u \,dD = 4\lambda \int_D \frac{vu}{(1+r^2)^2} \,dD\,.
\end{equation}

We use finite elements to approximate the eigenvalues and eigenfunctions of \eqref{eq:varP}. 
The general procedure we follow is:
\begin{enumerate}
\item Discretize the region  $D$ using a triangular mesh $\Tau_h=\bigcup_{i=1}^{N_h} \tau_i$, 
with a size of an individual triangle $\tau\in\Tau_h$ parameterized by $h>0$.
\item Introduce a finite-dimensional subspace $V_h$ of $\cal{H}$, consisting of finite element 
basis functions $\{\phi_i\}_{i=1}^{N_h}$ on $\Tau_h$;
\item Denote $(v_h$, $\lambda_h)\in V_h\times \RE$, with $v_h=(v_1, v_2, \dots,v_{N_h})^t$, the solution of the finite-dimensional generalized 
eigenvalue problem
\begin{equation}\label{eq:varPh}  
A_h{v}_h =\lambda_h B_h {v_h}\,, 
\end{equation}
where 
\begin{equation}
(A_h)_{ij}:= \int_{D} \nabla \phi_i \cdot \nabla \phi_j \,dD\,, 
\quad 
(B_h)_{ij} := \int_{D} \frac{4\phi_i  \phi_j}{(1+r^2)^2} \,dD\,. 
\end{equation}
Clearly, problem \eqref{eq:varPh} is the discrete analog of \eqref{eq:varP}. The eigenpair $(v_h, \lambda_h)$ is  computed using some iterative algorithm, until a prescribed tolerance $\mathtt{tol}$ is reached.
\item Steps 1-3 are repeated with smaller and smaller $h$ until some other stopping criterion is attained.
\end{enumerate}

We now present the results of some numerical experiments based on this strategy. 

\subsection{Conforming finite elements}
In the first set of experiments, the choice of approximating subspaces $V_h$ was a sequence of 
\emph{P1-conforming finite element spaces (Courant triangles)}, see \cite{Braess}. This means that for a given $h>0$, and a triangulation $\Tau_h$ of the domain,  
$$ 
V_h:= \{v\in \mathcal{H};\ v\vert_\tau = \text{ polynomial of degree } \leq 1 
\text{ for every } \tau \in \Tau_h\}\,.
$$ 
The discrete generalized eigenvalue problem for each $h$ was solved using an Arnoldi iteration with shift 2.6.
For details on the Arnoldi iteration, see, e.g., \cite{Golub, Trefethen}. 

\begin{figure}[!hbt]
\begin{center}
\includegraphics[width=0.9\textwidth,clip=true]{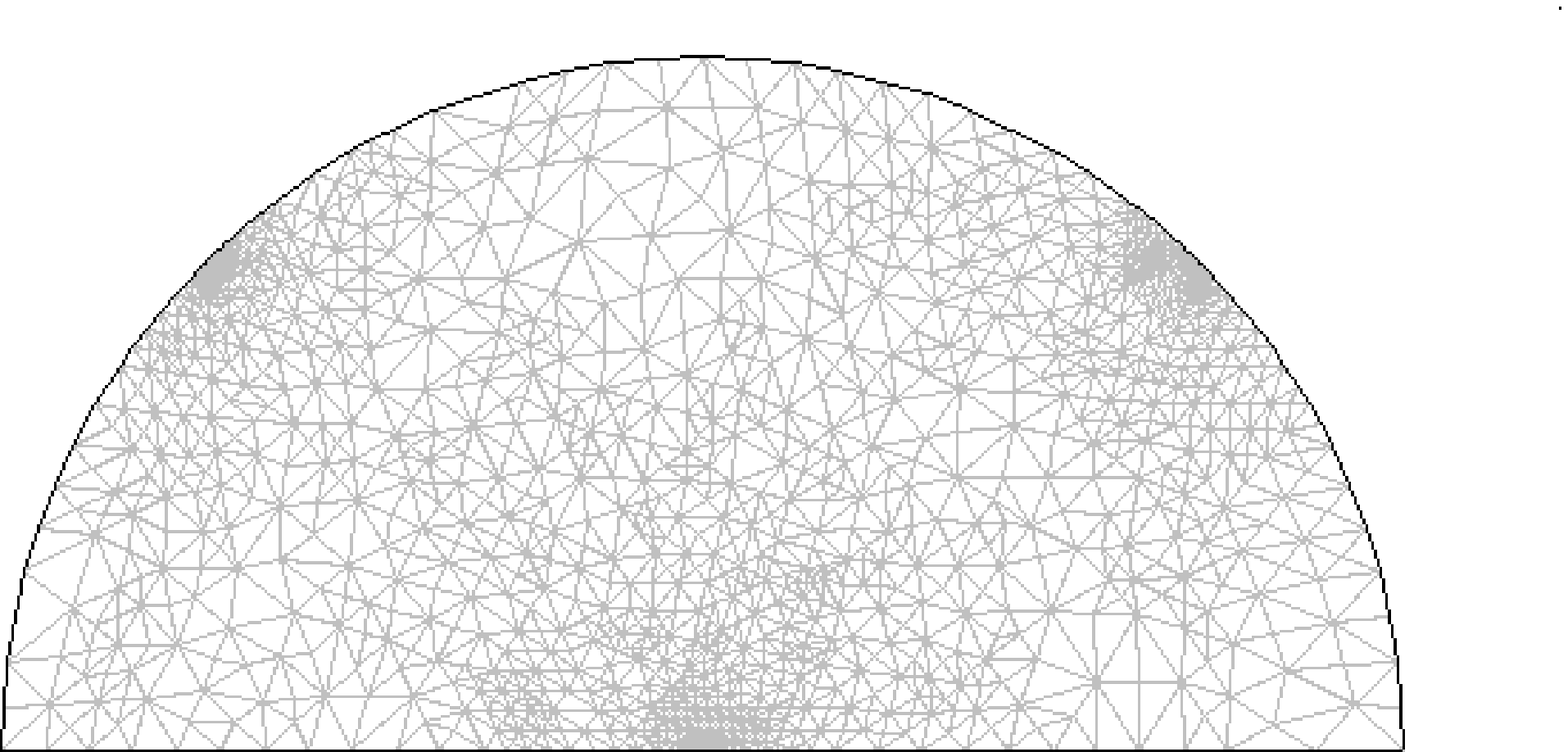}
\caption{A finite element mesh}
\label{fig:mesh}
\end{center}
\end{figure}

{\bf Experiment 1:}
The initial triangulation is based on a graded mesh, with more triangles located near the singularities of the 
eigenfunction (see Figure \ref{fig:mesh}).

The refinement strategy was based on simply subdividing each triangle in $\Tau_{h}$ into 4 while preserving 
the quality of the mesh, yielding a new triangulation $\Tau_{h/2}$.  The eigenvalue solver was run until a tolerance of $\mathtt{tol}=10^{-16}$ was achieved.  The computation was performed using {\tt FreeFem++} for generating the finite elements and the meshes, and {\tt ARPACK} for the eigenvalue solve.  The meshes were refined until the measure of error, 
$$ 
\mathtt{Err}:=\int_D |\nabla u_h|^2 - 4\lambda_h \frac {u_h^2}{(1+r^2)^2} \,dD\,,
$$ 
satisfied $|\mathtt{Err}|<5\times 10^{-10}$. The results are tabulated in Table \ref{Table:1}.

\begin{table} 
\begin{center}
 \begin{tabular}{|c|c|c|c|c|}  \hline
$\lambda_h$& Err& $N_h$ & No.  of& No. of \\ 
& & & Triangles & Arnoldi iterates\\ \hline
2.45590105457 & 2.00434573363e-05 & 169 & 288 & 11\\  \hline
2.36301118569 & 1.32592470511e-06 & 625 & 1152&11\\  \hline
2.32089716556 & 8.16135748742e-08 & 2401 & 4608&12\\  \hline
2.30111238184 & 4.80693483786e-09 & 9409 & 18432&11\\  \hline
2.29161462311 & 2.79565739833e-10 & 37249 & 73728&11\\  \hline
\end{tabular}
\caption{Using P1-conforming finite elements\label{Table:1}}
\end{center}
\end{table}

{\bf Experiment 2:}
This experiment was conducted using {\tt MATLAB}'s finite element package {\tt PDEToolbox}, and 
the eigenvalue solve was performed using {\tt ARPACK} routines. A sequence of triangular meshes 
was created, starting from the coarsest mesh, and refining 5 times. The major difference between 
this and the previous experiment is in the manner in which the zero Dirichlet data is enforced. 
 \begin{table}
 \begin{center}
\begin{tabular}{|c|c|c|c|c|} \hline
$\lambda_h$& $N_h$ & No.of \\ 
& &  Triangles \\ \hline
 2.55310562723060 & 77     &  126\\  \hline
  2.40400118356918    &   279       &  504\\  \hline
   2.33742285062686     &1061        &2016\\  \hline
     2.30582934149898     &4137        &8064\\  \hline
   2.29039772121374    &16337       &32256\\  \hline
       2.28276090970583   &64929      &129024\\  \hline
      2.27895954902635   &258881      &516096\\  \hline
\end{tabular}
\caption{Using P1-conforming finite elements in {\tt MATLAB}}
\end{center}
\end{table}

\subsection{Nonconforming finite elements}
In the second set of experiments, we used \emph{P1-nonconforming finite elements (Crouzeix elements)}, see \cite{Braess}.  These are 
defined as 
\begin{multline*}
V_h:=\{ v\in L^2(D); v|_\tau\text{ is linear for each }\tau \in \Tau_h,\\
v\text{ is continuous at the midpoints of triangle edges}\}
\end{multline*}
for a given $h>0$ and a triangulation $\Tau_h$. Note that $V_h$ is not a subspace of $\mathcal{H}$; for more information on 
the use of nonconforming elements in eigenvalue problems, see  \cite{Armentano}.
As before, the discrete generalized eigenvalue problem is solved using an Arnoldi 
iteration with 
a shift of 2.2 until a tolerance of $\mathtt{tol}=10^{-16}$ is achieved. The refinement strategy was to subdivide each triangle  in 
$\Tau_h$ into 4 subtriangles, yielding a new mesh $\Tau_{h/2}$. The meshes were refined until a measure of error  
$$ 
\mathtt{Err}:=\int_D |\nabla v_h|^2 - 4\lambda_h \frac {v_h^2}{(1+r^2)^2} \,dD\,,
$$ 
satisfied $|\mathtt{Err}|<5\times 10^{-10}$. 
The results are presented in Table \ref{Table:3}.
\begin{table}
\begin{center}
\begin{tabular}{|c|c|c|c|c|} \hline
$\lambda_h$& $\mathtt{Err}$& $N_h$ & No.  of& No. of \\ 
& & & Triangles & Arnoldi iterates\\ \hline
2.13042989031 & -1.5494060025e-05 & 169 & 288 &    11  \\
2.20743747322 & -7.89999667122e-07 & 625 & 1152 &  11     \\
2.24561396752 & -3.8455570927e-08 & 2401 & 4608 &   11    \\
2.26440630518 & -1.87263030138e-09 & 9409 & 18432 &   11    \\
2.27364314423 & -8.78059287464e-11 & 37249 & 73728 &  11     \\ \hline
\end{tabular}
\caption{Using P1-nonconforming finite elements\label{Table:3}}
\end{center}
\end{table}

\begin{figure}[!hbt]
\begin{center}
\includegraphics[width=0.9\textwidth,clip=true]{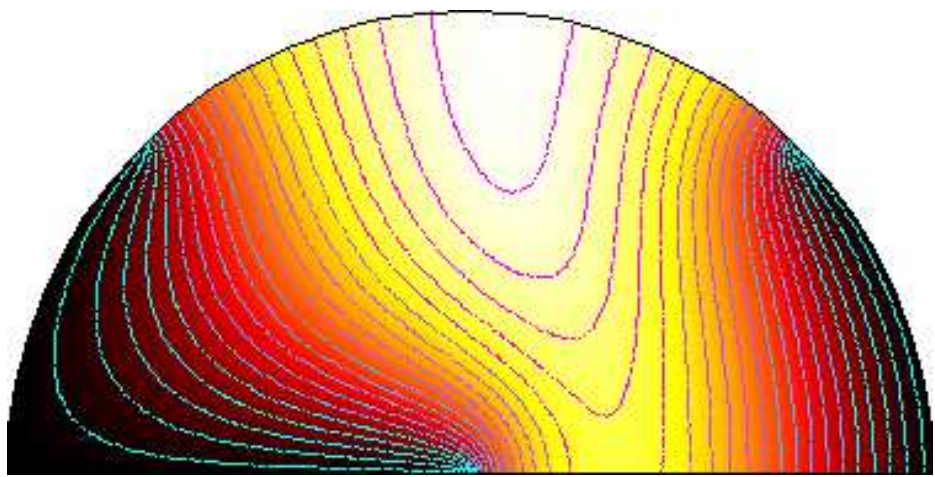}
\caption{Contour lines of the first eigenfunction}
\label{fig:contour}
\end{center}
\end{figure}

In each of the experiments above, we found that the computed eigenvalues appeared to converge to a value greater than 2.27. 
The associated eigenfunctions also appear to converge to a function whose contour lines are shown in Figure \ref{fig:contour}. 

\section*{Acknowledgements}   The authors would like to thank Herschel
Farkas, Eyal Goren, Daniel Grieser, Jacques Hurtubise, and Rafe Mazzeo for useful discussions and suggestions. 
Part of this work was completed when I.P. was visiting Heriot-Watt University and the Mathematical Research Institute at Oberwolfach,  M.L. was visiting Universit\'e de 
Montr\'eal, and D.J. was visiting Max Planck Institute for Mathematics in Bonn.  Hospitality of these institutions is greatly appreciated.  D.J., M.L. and I.P. also wish to thank the organizers of the LMS  Durham Symposium ``Operator Theory and Spectral Analysis'' for providing excellent conditions for collaboration.

The research of D.J. was partially supported by  NSERC, FQRNT, Dawson fellowship and Alfred P. Sloan Foundation fellowship.
The research of N.Nig. and I.P. was partially supported by NSERC and FQRNT.

\end{document}